\documentclass[sn-mathphys-num]{sn-jnl}
\usepackage{graphicx}
\usepackage{multirow}
\usepackage{amsmath,amssymb,amsfonts}
\usepackage{amsthm}
\usepackage{mathrsfs}
\usepackage[title]{appendix}
\usepackage{xcolor}
\usepackage{textcomp}
\usepackage{manyfoot}
\usepackage{booktabs}
\usepackage{algorithm}
\usepackage{algorithmicx}
\usepackage{algpseudocode}
\usepackage{listings}
\usepackage{hyperref}
\theoremstyle{thmstyleone}
\theoremstyle{thmstyletwo}
\theoremstyle{thmstylethree}

\newtheorem{exmp}{Example}[section]
\newtheorem{definition}{Definition}[section]
\newtheorem{theorem}{Theorem}[section]
\newtheorem{corollary}[theorem]{Corollary}

\newtheorem{lemma}[theorem]{Lemma}
\numberwithin{equation}{section}
\begin{document}
\title[Geometry of $\delta$-almost gradient Yamabe solitons on pseudo-Riemannian manifolds]{Geometry of $\delta$-almost gradient Yamabe solitons on pseudo-Riemannian manifolds}

\author[1]{\fnm{Rajdip} \sur{Biswas}}\email{rajdipbiswas467@gmail.com}

\author[2]{\fnm{Santu} \sur{Dey}}\email{santu.mathju@gmail.com}

\author[1]{\fnm{Arindam} \sur{Bhattacharyya}}\email
{bhattachar1968@yahoo.co.in}

\affil[1]{\orgdiv{Department of Mathematics}, \orgname{Jadavpur University}, \orgaddress{\city{Kolkata}, \postcode{700032}, \country{ India}}}

\affil[2]{\orgdiv{Department of Mathematics}, \orgname{Bidhan Chandra College}, \orgaddress{\city{Asansol}, \postcode{713304}, \country{India}}}

\abstract
    {In this article, we studied $\delta$-almost Yamabe solitons within the framework of paracontact metric manifolds. First, we proved that for a paracontact metric manifold $M$, if a paracontact metric $g$ represents a $\delta$-almost Yamabe soliton associated with the potential vector field $Z$ being an infinitesimal contact transformation, then $Z$ is Killing and if the potential vector field $Z$ is collinear with $\xi$, then the manifold $M$ is $K$-paracontact. 
     Next, if we take a $K$-paracontact metric manifold admitting $\delta$-almost Yamabe soliton with the potential vector field $Z$ parallel to the characteristic vector field and with constant scalar curvature then either scalar curvature will vanish or $g$ becomes a $\delta$-Yamabe soliton under a certain condition. 
     We established some results on $K$-paracontact manifold admitting $\delta$-almost gradient Yamabe soliton.
     Moreover, we consider a $(k,\mu)$-paracontact metric manifold admitting a non-trivial $\delta$-almost gradient Yamabe soliton. We shown that the potential vector field $Z$ is parallel to $\xi$. 
     We have also discussed about $\delta$-almost gradient Yamabe soliton on
     the para-Sasakian manifold. 
     Finally, we consider a para-cosymplectic manifold with a $\delta$-almost Yamabe soliton. 
     In the end, we construct two examples of $K$-paracontact metric manifolds with $\delta$-almost Yamabe soliton.}


\maketitle

\section{Introduction}
The significance of paracontact geometry is emitted from the para-Kähler manifold theory, while the theory of Legendre foliations is linked to the geometry of paracontact metric manifolds. Paracontact geometry methods play an important role in modern mathematics. The geometry of almost paracontact manifolds is a natural extension of the almost para Hermitian geometry, just like almost contact manifolds correspond to the almost Hermitian ones. Last some years, many authors studied on paracontact geometry has evolved from the mathematical formalism of classical mechanics (see \cite{ghosh}). On the analogy of almost contact manifolds, Sato \cite{sato} introduced the notion of the almost paracontact manifolds. An almost contact manifold is always odd dimensional, but an almost paracontact manifold could be of even dimensional as well. Takahashi \cite{taka} defined almost contact manifolds, in particular, Sasakian manifolds equipped with an associated pseudo-Riemannian metric. Later, Kaneyuki and Williams \cite{kan} introduced the notion of an almost paracontact pseudo-Riemannian structure, as a natural odd dimensional counterpart to paraHermitian structure. In \cite{zam}, Zamkovoy showed that any almost paracontact structure admits a pseudo-Riemannian metric with signature $(n+1,n)$. In recent years, almost paracontact structure has been studied by many authors, particularly since the appearance of \cite{zam}. A systematic study of paracontact metric manifolds had been started in the paper \cite{zam}. The technical apparatus introduced in that paper is essential for further investigations about paracontact metric geometry.\par
The concept of Yamabe flow was first introduced by Hamilton \cite{Ham1} to construct Yamabe metrics on compact Riemannian manifolds. On a Riemannian or
pseudo-Riemannian manifold $M$, a time-dependent metric $g(\cdot,t)$ is said to evolve
by the Yamabe flow if the metric $g$ satisfies the given equation,
\begin{equation}
 \frac{\partial g(t)}{\partial t}=-r(g(t)), \hspace{5pt} g(0)=g_{0},\label{1.1}
\end{equation}
where $r$ is the scalar curvature of the manifold $M$.\\
In 2-dimension the Yamabe flow is equivalent to the Ricci flow \cite{Ham2} (defined by $\frac{\partial g(t)}{\partial t} =-2S(g(t))$, where S denotes the Ricci tensor). But in dimension $\geq$3
the Yamabe and Ricci flows do not agree, since the Yamabe flow preserves the conformal class of the metric but the Ricci flow does not in general. A Yamabe soliton \cite{Bar} correspond to self-similar solution of the Yamabe flow, is defined on a
Riemannian or pseudo-Riemannian manifold $(M, g)$ as,
\begin{equation}
  \frac{1}{2}\mathfrak{L}_{Z}g=(r-\lambda)g,\label{1.2}
\end{equation}
where $\mathfrak{L}_{Z}g$ denotes the Lie derivative of the metric $g$ along the vector field $Z$, $r$ is the scalar curvature and $\lambda$ is a constant. Moreover Yamabe soliton is said to be expanding, steady, shrinking depending on $\lambda$ being positive, zero, negative 
respectively. If $\lambda$ is a smooth function then \eqref{1.2} is called almost Yamabe soliton
\cite{Bar}.\\
Very recently, Chen \cite{Che} introduced a new concept, named $\delta$-almost Yamabe soliton. According to Chen, a Riemannian metric is said to be a $\delta$-almost Yamabe soliton if there exists a smooth vector field $Z$, a $C^{\infty}$ function $\lambda$ and a nonzero function $\delta$ such that,
\begin{equation}
    \frac{\delta}{2}\mathfrak{L}_{Z}g = (r-\lambda)g \label{1.3}
\end{equation}
holds. If for any smooth function $u$, $Z=\nabla u$ then the previous equation is called $\delta$-almost gradient Yamabe soliton. If $\lambda$ is constant, then \eqref{1.3} is called $\delta$-Yamabe soliton equation. For more information on Yamabe solitons, gradient Yamabe solitons and their generalizations within the framework of paracontact Riemannian geometry, we recommend the papers \cite{Ro, Sarkar, Kumara, Li, KU, De, Bla} and the references cited there. \par
A fascinating problem in differential geometry is to obtain some conditions on a Riemannian manifold endowed with a soliton structure so that it is isometric to an Euclidean sphere. For several sufficient conditions under which an almost Ricci soliton on Riemannian or contact metric manifolds are isometric to an Euclidean sphere the readers can see \cite{bar, desh, desh1}. For instance, Ricci soliton and gradient Ricci solitons have been demonstrated by many authors on paracontact metric structure (see \cite{Bejan, Calvaruso, iva, liu, martin, perrone2016some}). In particular, Calvaruso et al. \cite{Calvaruso} explicitly exhibited Ricci solitons on 3-dimensional almost paracontact manifolds. It is a natural odd-dimensional counterpart to para-Hermitian structures, just like contact metric structures correspond to almost Kähler ones. Ricci solitons and their generalizations have been enormously studied by some authors within the framework of contact and paracontact metric manifolds (see in details \cite{roy*, dey*, SD, SD1, SD2, SD3, SD4, SD5}). Recently, I. K. Erken \cite{Erken} demonstrated Yamabe solitons on 3-dimensional para-cosymplectic manifold and proved some vital results like the manifold is either $\eta$-Einstein or Ricci flat. Moreover, we can get more latest concept equipped with soliton geometry in \cite{suh1, suh2, suh3, suh4, Li1, Li2, siraj}. Motivated by the above studies, we consider the $\delta$-almost Yamabe soliton and $\delta$-almost gradient Yamabe soliton on the pseudo-Riemannian manifold. \\
The organization of this article is as follows: Sections 1 and 2 focus on introducing and explaining the fundamental concepts of paracontact metric($K$-paracontact and para-Sasakian) manifolds, which consists of basic definitions and notions. In section 3, we investigate $\delta$-almost Yamabe soliton and $\delta$-almost gradient Yamabe soliton on $K$-paracontact, para-Sasakian and $(k,\mu)$ paracontact metric manifolds in connection several conditions of potential vector field. Section 4 provides some characterizations of $\delta$-almost Yamabe soliton on para-cosymplectic manifold and classify a result on scalar curvature. In section 5, we give some examples to support our findings.

\section{Preliminaries}
This section provides a brief overview of the concepts which were used throughout the paper.\\
A smooth manifold $M^{n}$ has an almost paracontact structure $(\phi,\xi,\eta)$ if it is endowed with a vector field $\xi$ (called the Reeb vector field), a (1,1)-type tensor field $\phi$, and a 1-form $\eta$ obeying the subsequent conditions:\\
\begin{equation}
  \phi ^2 = I - \eta \otimes \xi , \label{2.1}
  \end{equation}
  \begin{equation}
    \eta (\xi ) = 1,\label{2.2}
  \end{equation}
  \begin{equation}
  \hspace{4pt}\phi(\xi)=0,\hspace{4pt}\eta o \phi=0 .\label{2.3}
\end{equation}
If an almost paracontact manifold has a pseudo-Reimannian metric $g$ in which
\begin{equation}
    g(\phi X_1,\phi X_2)=-g(X_1,X_2) + \eta(X_1)\eta(X_2), \hspace{5pt}X_1,X_2\in\chi(M),\label{2.4}
\end{equation}
where $\chi(M)$ is the module over $C^{\infty}(M)$ of all vector fields on $M$, then $M$ admits an almost paracontact metric structure $(\phi,\xi,\eta,g)$ and $g$ is referred to as a compatible metric.
The Nijenhuis torsion is defined by
             \begin{equation}
             [\phi,\phi](X_1,X_2)=[\phi X_1,\phi X_2]+\phi^2[X_1,X_2]-\phi[X_1,\phi X_2]-\phi[\phi X_1,X_2] \label{2.5}
             \end{equation}
        for all $X_1,X_2\in \chi(M)$.
The almost paracontact metric manifold is called normal if $N_{\phi}=[\phi,\phi]-2d\eta\otimes\xi$ vanishes. 
The \textit{fundamental 2-form} $\Phi$ of an almost paracontact metric structure $(\phi,\xi,\eta,g)$ is given by $\Phi(X_1,X_2)=g(X_1,\phi X_2)$ for any vector fields $X_1$ and $X_2$ on $M$. In case that $\Phi=d\eta$, then $M^{n}(\phi,\xi,\eta,g)$ is referred to as a paracontact metric manifold.
Here we can define two self-adjoint operators on $M^{n}(\phi,\xi,\eta,g)$ as $l=R(.,\xi)\xi$ and $h=\frac{1}{2}\mathfrak{L}_{\xi}\phi$, where $\mathfrak{L}_{\xi}$ is the \textit{Lie derivative} in the direction of $\xi$ and $R$ is the \textit{Riemann curvature tensor} of $g$ defined as follows:
\begin{equation}
    R(X_1,X_2)=[\nabla_{X_1},\nabla_{X_2}]-\nabla_{[X_1,X_2]},\hspace{5pt} X_1,X_2\in\chi(M),\label{2.6}
\end{equation}
where $\nabla$ is the operator of covariant differentiation of $g$.\\
On a paracontact metric manifold, the following formulae hold \cite{zam}:
\begin{align}
    \nabla_{X_1}\xi &= -\phi X_1 + \phi hX_1, \hspace{5pt} X_1\in
    \chi(M),\label{2.7}\\
     \mathcal{S}(\xi,\xi) &= g(Q\xi,\xi)=Tr_{g}l=Tr_{g}(h^2)-2n, \label{2.8}   
\end{align}
where $Q$ stands for the Ricci operator associated to the Ricci tensor defined as $\mathcal{S}(X_1,X_2)=g(QX_1,X_2)$ for $X_1,X_2\in\chi(M)$. If $h=0$, i.e., the vector field $\xi$ is Killing, then $M$ is referred to as a \textit{K-paracontact manifold}. On \textit{K-paracontact manifold}, the following formulae hold \cite{zam}:
\begin{align}
\nabla_{X_1}\xi&=-\phi X_1,\label{2.9}\\
 R(X_1,\xi)\xi&=-X_1+\eta(X_1)\xi,\label{2.10}\\
 Q\xi&=-2n\xi\label{2.11}
\end{align}
for all $X_1\in\chi(M)$.
\begin{definition}
    If the normality condition is satisfied in a paracontact metric manifold, i.e. satisfies $[\phi, \phi] - 2d\eta\otimes\xi = 0$ then it is called a para-Sasakian manifold. This condition is equivalent to
    \begin{equation}
        (\nabla_{X_1}\phi)X_2 = -g(X_1,X_2)\xi +\eta(X_2)X_1.\label{2.12}
    \end{equation}
\end{definition}
    A para-Sasakian manifold is in particular \textit{$K$-paracontact}. The converse holds in dimension 3 but not in general for higher dimensions. Every para-Sasakian manifold satisfies
    \begin{align}
        R(X_1,X_2)\xi &=\eta(X_1)X_2-\eta(X_2)X_1, \label{2.13}\\
        R(X_1, \xi)X_2 &=g(X_1,X_2)\xi-\eta(X_2)X_1.\label{2.14}
    \end{align}

\begin{definition}
   \cite{Cappe} If the curvature tensor $R$ in a paracontact metric manifold obeys
         \begin{equation}
        R( X_1,X_2)\xi=k(\eta(X_2)X_1-\eta(X_1)X_2)+\mu(\eta(X_2)hX_1-\eta(X_1)hX_2),\label{2.15}
         \end{equation}
         then it is called a $(k,\mu)$-paracontact metric manifold, where $k,\mu$ are real and $2h$ is the Lie differentiation of $\phi$ in the direction of $\xi$.
\end{definition}

\begin{definition}
A vector field $Z$ on a contact manifold is called an infinitesimal contact transformation if it preserves the contact form $\eta$ meaning that
\begin{equation}\label{2.16}
    \mathfrak{L}_{Z}\eta=\sigma\eta
\end{equation}
for some smooth function $\sigma$ on $M$. If $\sigma=0$, then $Z$ is referred to as a strictly infinitesimal contact transformation.
\end{definition}
Let's denote the volume element $\eta\wedge(d\eta)^{m}$ as $\Omega$, where $n=2m$ or $2m-1$. According to the definition of a paracontact structure, we have $\Omega\neq 0$. Now, applying exterior differentiation to equation \eqref{2.16} we obtain
\begin{equation}\label{2.17}
\mathfrak{L}_{Z}d\eta=d\mathfrak{L}_{Z}\eta=d\sigma\wedge\eta+\sigma d\eta.
\end{equation}
Since $\Omega$ is a volume form, $\Omega$ is closed. Thus, using Cartan’s formula we find
\begin{equation}\label{2.18}
    \mathfrak{L}_{Z}\Omega=(div Z)\Omega.
\end{equation}
So, we take the Lie derivative of $\Omega=\eta\wedge(d\eta)^{m}$ 
and using equations \eqref{2.17} and \eqref{2.18} deduce that
\begin{equation}\label{2.19}
    divZ=(m+1)\sigma.
\end{equation}
These equations are useful in the next sections.

\section{Proof of the main outcomes}
In this section, we will discuss our results and provide their proofs.
Before introducing the first theorem in this section, we will establish the following lemma: \vspace{6pt}

\begin{lemma}
 If the metric $g$ of a paracontact metric manifold represents a $\delta$-almost Yamabe soliton, then the following properties hold.
    \begin{align}
        \eta(\mathfrak{L}_{Z}\xi) &= \frac{\lambda -r}{\delta},\label{3.37}\\
        (\mathfrak{L}_{Z}\eta)\xi &= \frac{r-\lambda}{\delta}.\label{3.38}
    \end{align}
\end{lemma}
\proof
From the soliton equation \eqref{1.3} we get, 
$\delta(\mathfrak{L}_{Z}g)(X_1,\xi)=2(r-\lambda)\eta(X_1)$.Using of this in the Lie-derivative of $\eta(X_1)=g(X_1,\xi)$ along $Z$ we get
\begin{equation}\label{3.39}
(\mathfrak{L}_{Z}\eta)(X_1)-g(X_1,\mathfrak{L}_{Z}\xi)=2(\frac{r-\lambda}{\delta})\eta(X_1).
\end{equation}
Now, we operate $\mathfrak{L}_{Z}$ to $g(\xi,\xi)=1$ and get
\begin{align}
    \eta(\mathfrak{L}_{Z}\xi)&=-(\frac{r-\lambda}{\delta}).\nonumber
\end{align}
Using this with equation \eqref{3.39} we obtain
\begin{equation}
    (\mathfrak{L}_{Z}\eta)\xi=(\frac{r-\lambda}{\delta}).\nonumber
\end{equation}
This completes the proof.\\

Now we use the above lemma in the next theorem.
\begin{theorem}
Let $M^n$ be a paracontact metric manifold. If g is a $\delta$-almost Yamabe soliton with the potential vector field $Z$ being an infinitesimal contact transformation, then $Z$ is Killing.
\end{theorem}
\proof From \eqref{3.39} and \eqref{2.16}, we find that
\begin{equation}
    g(X_1,\mathfrak{L}_{Z}\xi)=\sigma\eta(X_1)-2(\frac{r-\lambda}{\delta})\eta(\xi).
\end{equation}
By substituting $X_1$ with $\xi$ in this equation and using  \eqref{3.37}, we get 
$\sigma=\frac{r-\lambda}{\delta}$. From the equation \eqref{1.3}, we can find that $divZ=n\frac{(r-\lambda)}{\delta}$.
Using this result in equation \eqref{2.19} confirms that $r=\lambda$. Thus, the proof is complete.

\begin{theorem}
    If pseudo-Riemannian metric $g$ of a paracontact metric manifold $M^n$ is a $\delta$-almost Yamabe soliton with potential vector field $Z$ being parallel to Reeb vector field, then $M^n$ is K-paracontact.
\end{theorem}
\proof Let $Z=f\xi$, for non-zero function $f$ on $M$. Using this in \eqref{1.3} and \eqref{2.7} we get
\begin{equation}\label{3.41}
    \delta[(X_{1}f)\eta(X_2)+(X_{2}f)\eta(X_1)+2fg(\phi hX_1,X_2)]=2(r-\lambda)g(X_1,X_2).
\end{equation}
Now by substituting $X_1$ by $\phi X_1$ and $X_2$ by $\xi$ in \eqref{3.41} we get
\begin{equation}
    \phi\nabla f=0.\nonumber
\end{equation}
By operating $\phi$ on this equation shows that
\begin{equation}
    \nabla f=(\xi f)\xi.\nonumber
\end{equation}
Again putting $X_1=X_2=\xi$ in \eqref{3.41} we derive
\begin{align}
    \xi f &=\frac{r-\lambda}{\delta},\nonumber\\
    \nabla f &=(\frac{r-\lambda}{\delta})\xi.\nonumber
\end{align}
Now differentiating this along $X_1$ and using \eqref{2.7} we obtain
\begin{equation}
    \nabla_{X_1}\nabla f = X_1(\frac{r-\lambda}{\delta})\xi-(\frac{r-\lambda}{\delta})\phi X_1 + (\frac{r-\lambda}{\delta})\phi hX_1.\nonumber
\end{equation}
We know that, $g(\nabla_{X_1}\nabla f,X_2)=g(\nabla_{X_2}\nabla f,X_1)$.
Therefore we get
\begin{equation}
    X_1(\frac{r-\lambda}{\delta})\eta(X_2)-X_2(\frac{r-\lambda}{\delta})\eta(X_1)-2(\frac{r-\lambda}{\delta})g(\phi X_1,X_2)=0. \nonumber
\end{equation}
Putting $X_1=\phi X_1$ and $X_2=\phi X_2$ in the last equation, we achieve 
\begin{equation}
    r=\lambda. \nonumber
\end{equation}
Therefore it implies that $f$ is constant. Hence from \eqref{1.3} we see that $\xi$ is a Killing vector. This completes the proof. 
\\

Now we focus on $K$-paracontact manifold admitting $\delta$-almost Yamabe soliton and $\delta$-almost gradient Yamabe soliton, we have the following result:
\begin{theorem}
     Assume that $M^{n}(\phi,\xi,\eta,g)$ is a $K$-paracontact metric manifold admitting a $\delta$-almost Yamabe soliton with the potential vector field $Z$ parallel to the Reeb vector field $\xi$ and with constant scalar curvature. Then either scalar curvature will vanish or $g$ is a $\delta$-Yamabe soliton, whenever $\lambda=\delta$.
\end{theorem} 
\proof
Taking covariant derivative of \eqref{1.3} along $X_{3}\in\chi(M)$ 
\begin{equation}\label{3.1}
    (X_{3}\delta)\mathfrak{L}_{Z}g(X_1,X_2)+\delta(\nabla_{X_3}\mathfrak{L}_{Z}g)(X_1,X_2)=2(-X_3\lambda)g(X_1,X_2) .
\end{equation}
On pseudo-Reimannian manifold $(M,g)$ we have
\begin{equation}\label{3.2}
    (\nabla_{X_3}\mathfrak{L}_{Z}g)(X_1,X_2)=g((\mathfrak{L}_{Z}\nabla)(X_3,X_1),X_2)+g((\mathfrak{L}_{Z}\nabla)(X_3,X_2),X_1).
\end{equation}
Since $g$ is parallel. By utilizing the symmetry of the (1,2)-type tensor field, $\mathfrak{L}_{Z}\nabla$, i.e;
$\mathfrak{L}_{Z}\nabla(X_1,X_2)=\mathfrak{L}_{Z}\nabla(X_2,X_1)$ and interchanging the roles of $X_1,X_2,X_3$ in the preceding equation, we can compute
\begin{equation}\label{3.3}
\begin{split}
    g((\mathfrak{L}_{Z}\nabla)(X_1,X_2),X_3)=\frac{1}{2}(\nabla_{X_1}\mathfrak{L}_{Z}g)(X_2,X_3)+\frac{1}{2}(\nabla_{X_2}\mathfrak{L}_{Z}g)(X_3,X_1)\\-\frac{1}{2}(\nabla_{X_3}\mathfrak{L}_{Z}g)(X_1,X_2) .
\end{split}
\end{equation}
Next by combining \eqref{3.1} and \eqref{3.3} we deduce

\begin{equation}\label{3.4}
    \begin{split}
     \delta g((\mathfrak{L}_{Z}\nabla)(X_1,X_2),X_3)=(X_3\lambda)g(X_1,X_2)-(X_1\lambda)g(X_2,X_3)-(X_2\lambda)g(X_3,X_1)\\
     -\frac{(r-\lambda)}{\delta}[(X_1\delta)g(X_2,X_3)+
     (X_2\delta)g(X_3,X_1)-(X_3\delta)g(X_1,X_2)].
    \end{split}
\end{equation}
Now by substituting $\xi$ for $X_2$ in \eqref{3.4} we obtain
\begin{equation}\label{3.5}
\begin{split}
    \delta g((\mathfrak{L}_{Z}\nabla)(X_1,\xi),X_3)=(X_3\lambda)\eta(X_1)-(X_1\lambda)\eta(X_3)-(\xi\lambda)g(X_3,X_1)\\-\frac{(r-\lambda)}{\delta}[(X_1\delta)\eta(X_3)+(\xi\delta)g(X_3,X_1)-(X_3\delta)\eta(X_1)].
\end{split}
\end{equation}
Setting $X_2=\xi$ in the well-known formula 
\begin{equation}\label{3.6}
    \nabla_{X_1}\nabla_{X_2}Z-\nabla_{\nabla_{X_1}X_2}Z-R(Z,X_1)X_2=(\mathfrak{L}_{Z}\nabla)(X_1,X_2)
\end{equation}
and by plugging the value of $(\mathfrak{L}_{Z}\nabla)(X_1,\xi)$ in \eqref{3.5} we obtain
\begin{equation}\label{3.7}
\begin{split}
\delta g(\nabla_{X_1}\nabla_{\xi}Z-\nabla_{\nabla_{X_1}\xi}Z-R(Z,X_1)\xi,X_3)=(X_3\lambda)\eta(X_1)-(X_1\lambda)\eta(X_3)-\\(\xi\lambda)g(X_3,X_1)-\frac{(r-\lambda)}{\delta}[(X_1\delta)\eta(X_3)+(\xi\delta)g(X_3,X_1)-(X_3\delta)\eta(X_1)].
\end{split}
\end{equation}
Next taking $\xi$ instead of $X_1$ and $X_3$ we get
\begin{equation}\label{3.8}
    \delta g(\nabla_{\xi}\nabla_{\xi}Z-\nabla_{\nabla_{\xi}\xi}Z-R(Z,\xi)\xi,\xi)=-(\xi\lambda)-\frac{(r-\lambda)}{\delta}(\xi\delta).
\end{equation}
Assume the potential vector field $Z$ is a Jacobi field in direction of $\xi$, i.e., $\nabla_{\xi}\nabla_{\xi}Z-R(Z,\xi)\xi=0$. It follows from \eqref{3.8} that
\begin{equation}\label{3.9}
    (\xi\lambda)+\frac{(r-\lambda)}{\delta}(\xi\delta)=0.
\end{equation}
Now If $\lambda$ and $\delta$ both are equal then we get from \eqref{3.9} that
    
\begin{equation*}
     (1+\frac{r-\lambda}{\lambda})(\xi\lambda)=0,
\end{equation*}
\begin{equation}\label{3.10}
        r(\xi\lambda)=0.
\end{equation}
Therefore it follows from \eqref{3.10} that, either $r=0$ or $\lambda$ is constant which implies that $g$ becomes a $\delta$-Yamabe soliton. 

\begin{theorem}
    Suppose that $M^{n}(\phi,\xi,\eta,g)$ is a K-paracontact manifold admitting a  $\delta$-almost Yamabe soliton with the non-zero potential vector field $Z$ is collinear with $\xi$. Then the potential vector field $Z$ is a constant multiple of $\xi$.
\end{theorem}
 
\proof We know that
\begin{equation}\label{3.11}
    (\mathfrak{L}_{Z}g)(X_1,X_2)=g(\nabla_{X_1}Z,X_2)+g(X_1,\nabla_{X_2}Z).
\end{equation}
Since the potential vector field $Z$ is collinear with $\xi$, i.e., $Z=\sigma\xi$ for a non-zero differentiable function $\sigma$ on $M$.
Using the equation \eqref{2.9} in \eqref{3.11} we can compute
\begin{equation}\label{3.12}
    (\mathfrak{L}_{Z}g)(X_1,X_2)=X_1(\sigma)\eta(X_2)+X_2(\sigma)\eta(X_1).
\end{equation}
From equation \eqref{1.3} and \eqref{3.12} we get
\begin{equation}\label{3.13}
    \delta[X_1(\sigma)\eta(X_2)+X_2(\sigma)\eta(X_1)]=2(r-\lambda)g(X_1,X_2) .
\end{equation}
Putting $X_2=\xi$ in \eqref{3.13}
\begin{equation}\label{3.14}
    \delta X_1(\sigma)=[2(r-\lambda)-\delta\xi(\sigma)]\eta(X_1).
\end{equation}
Again putting $X_1=X_2=\xi$ in \eqref{3.13} we obtain
\begin{equation}
    \delta\xi(\sigma)=(r-\lambda).\label{3.15}
\end{equation}
Now from \eqref{3.14} and \eqref{3.15} we get
\begin{align}
    X_1(\sigma)&=\xi(\sigma)\eta(X_1),\nonumber\\
    g(\nabla\sigma,X_1)&=\xi(\sigma)g(X_1,\xi),\nonumber\\
    \nabla\sigma&=\xi(\sigma)\xi.\label{3.16}
\end{align}
From the equation \eqref{3.16} we can write $d\sigma=\xi(\sigma)\eta$, where $d$ denote the exterior differentiation. Taking the exterior derivative of the previous equation yields $d^2\sigma=d(\xi(\sigma))\wedge\eta+(\xi(\sigma))d\eta$. Using \textit{Poincare lemma} $d^2\equiv 0$ and then executing wedge product with $\eta$ in the foregoing equation we obtain
\begin{equation}\label{3.17}
    (\xi(\sigma))\eta\wedge d\eta=0.
\end{equation}
Hence $\xi(\sigma)=0$, since $\eta\wedge d\eta\neq 0$ in $M^{n}$. Thus we get $d\sigma=0$ which implies that $\sigma$ is constant.  
Therefore $Z$ is a constant multiple of $\xi$. This completes the proof.

\begin{lemma}
    Let $M^n$ be a paracontact metric manifold. If $g$ is a $\delta$-almost gradient Yamabe soliton, then the following properties hold:
    \begin{align}
        R(X_1,X_2)\nabla u &= X_1(\frac{r-\lambda}{\delta})X_2-X_2(\frac{r-\lambda}{\delta})X_1, \label{3.42}\\
        \mathcal{S}(X_2,\nabla u) &=-X_2(\frac{r-\lambda}{\delta})2n .\label{3.43}
    \end{align} 
\end{lemma}
\proof Since $g$ is $\delta$-almost gradient Yamabe soliton, so from \eqref{1.3} we obtain
\begin{equation}
    \delta \nabla_{X_1}\nabla u=(r-\lambda)X_1.\label{3.45}
\end{equation}
Now, from the expression of curvature tensor $R(X_1,X_2)=[\nabla_{X_1},\nabla_{X_2}]-\nabla_{[X_1,X_2]}$ we get
\begin{equation}
    R(X_1,X_2)\nabla u=X_1(\frac{r-\lambda}{\delta})X_2-X_2(\frac{r-\lambda}{\delta})X_1. \nonumber
\end{equation}
Again we can easily obtain \eqref{3.43} by contracting \eqref{3.42} over $X_1$.

\begin{theorem}
    If a K-paracontact manifold $M^n$ admits a $\delta$-almost gradient Yamabe  soliton with potential function $u$, then $\frac{r-\lambda}{\delta}-u=constant$.
\end{theorem}
\proof Firstly taking scalar product of \eqref{3.42} with $\xi$ we get
\begin{equation}
    g(R(X_1,X_2)\nabla u,\xi)=X_1(\frac{r-\lambda}{\delta})\eta(X_2)-X_2(\frac{r-\lambda}{\delta})\eta(X_1).\label{3.46}
\end{equation}
Now, taking inner product of \eqref{2.10} with $\nabla u$ yields
\begin{equation}
    g(R(X_1,\xi)\nabla u,\xi)=(X_{1}u)-\eta(X_1)(\xi u). \nonumber
\end{equation}
Hence from the above equation with \eqref{3.46} we obtain
\begin{equation}
    X_1(\frac{r-\lambda}{\delta})-\xi(\frac{r-\lambda}{\delta})\eta(X_1)=(X_{1}u)-\eta(X_1)(\xi u).\nonumber
\end{equation}
We can write this as $d(\frac{r-\lambda}{\delta}-u)=\xi(\frac{r-\lambda}{\delta}-u)\eta$. Taking the exterior derivative of the previous equation and using Poincare lemma yields $d(\xi(\frac{r-\lambda}{\delta}-u))\wedge \eta+\xi(\frac{r-\lambda}{\delta}-u)d\eta=0$. Executing wedge product with $\eta$ in the foregoing equation we obtain $\xi(\frac{r-\lambda}{\delta}-u)\eta \wedge d\eta=0$. Since $\eta\wedge d\eta\neq 0 $, we obtain $d(\frac{r-\lambda}{\delta}-u)=0$. It implies that, $\frac{r-\lambda}{\delta}-u=constant$.
\vspace{10pt}

Next consider $(k,\mu)$-paracontact and para-Sasakian manifold with $\delta$-almost gradient Yamabe soliton and prove the following conclusions.
\begin{theorem}
     Let $M^{n}$ be a $(k,\mu)$-paracontact metric manifold with $(k+\mu h)\neq 0$ whose metric is a non-trivial $\delta$-almost gradient Yamabe soliton. Then the potential vector field $Z$ is parallel to $\xi$.
\end{theorem}
\proof

Applying the covariant derivative of \eqref{3.45} along $X_2$, we obtain
\begin{equation}\label{3.29}
    \delta\nabla_{X_2}\nabla_{X_1}\nabla u=(X_2(r-\lambda))X_1+(r-\lambda)\nabla_{X_2}X_1-\frac{1}{\delta}(X_2\delta)(r-\lambda)X_1 ,
\end{equation}
where $X_1,X_2$ are any smooth vector fields.
Exchanging $X_1$ and $X_2$ in \eqref{3.29} we get
\begin{equation}\label{3.30}
    \delta\nabla_{X_1}\nabla_{X_2}\nabla u=(X_1(r-\lambda))X_2+(r-\lambda)\nabla_{X_1}X_2-\frac{1}{\delta}(X_1\delta)(r-\lambda)X_2 
\end{equation}
and
\begin{equation}\label{3.31}
    \delta\nabla_{[X_1,X_2]}\nabla u=(r-\lambda)[X_1,X_2].
\end{equation}
Now applying the equations \eqref{3.29}, \eqref{3.30} and \eqref{3.31} in the expression $R(X_1,X_2)X_3=\nabla_{X_1}\nabla_{X_2}X_3-\nabla_{X_2}\nabla_{X_1}X_3-\nabla_{[X_1,X_2]}X_3$, $X_3\in\chi(M)$ we obtain
\begin{equation}\label{3.32}
    \delta^{2}R(X_1,X_2)\nabla u=\delta[(X_1(r-\lambda))X_2]-\delta[(X_2(r-\lambda))X_1]-(X_1\delta)(r-\lambda)X_2+(X_2\delta)(r-\lambda)X_1.
\end{equation}
Executing the inner product of \eqref{3.32} with $\xi$ yields
\begin{equation}\label{3.33}
\begin{split}
     \delta^{2}g(R(X_1,X_2)\nabla u,\xi)=\delta[(X_1(r-\lambda))\eta(X_2)]-\delta[(X_2(r-\lambda))\eta(X_1)]-\\(X_1\delta)(r-\lambda)\eta(X_2)+(X_2\delta)(r-\lambda)\eta(X_1).
\end{split}
\end{equation}
Again we obtain from the relation \eqref{2.15} that
\begin{equation}\label{3.34}
\begin{split}
 \delta^{2}g(R(X_1,X_2)\nabla u, \xi)=-\delta^{2}k[X_1(u)\eta(X_2)-X_2(u)\eta(X_1)]-\delta^{2}\mu[(hX_1)(u)\eta(X_2)\\-(hX_2)(u)\eta(X_1)] . 
 \end{split}
\end{equation}
Now from equations \eqref{3.33} and \eqref{3.34} we get
\begin{equation}\label{3.35}
\begin{split}
  \delta[(X_1(r-\lambda))\eta(X_2)]-\delta[(X_2(r-\lambda))\eta(X_1)]-(X_1\delta)(r-\lambda)\eta(X_2)+(X_2\delta)(r-\lambda)\eta(X_1)\\=-\delta^{2}k[X_1(u)\eta(X_2)-X_2(u)\eta(X_1)]-\delta^{2}\mu[(hX_1)(u)\eta(X_2)-(hX_2)(u)\eta(X_1)].
  \end{split}
\end{equation}
Putting $X_2=\xi$ in \eqref{3.35} we get
\begin{equation}\label{3.36}
   Z=\nabla u=\frac{k}{k+\mu h} \xi(u)\xi.
\end{equation}
Therefore, it follows that if $(k+\mu h)\neq 0$ then $Z$ is parallel to $\xi$. This is the required proof.

\begin{theorem}
     Let the para-Sasakian manifold $M^{n}$ admits a nontrivial $\delta$-almost gradient Yamabe soliton. Then, the potential function of the $\delta$-almost gradient Yamabe soliton is constant.
\end{theorem}

\proof 
Let a para-Sasakian manifold $M^{n}$ admits a $\delta$-almost gradient Yamabe soliton $(g,u,\delta,\lambda)$. \\
In this case, the proof of the first part of this theorem is similar to the previous one. 
Again from the relation \eqref{2.13}, for para-Sasakian manifold we obtain
\begin{equation}\label{3.24}
    \delta^{2}g(R(X_1,X_2)\nabla u,\xi)=\delta^{2}[(X_{1}u)\eta(X_2)-(X_{2}u)\eta(X_1)].
\end{equation}
Combining equations \eqref{3.33} and \eqref{3.24} we get\\
\begin{equation}\label{3.25}
\begin{split}
    \delta[(X_1(r-\lambda))\eta(X_2)]-\delta[(X_2(r-\lambda))\eta(X_1)]-(X_1\delta)(r-\lambda)\eta(X_2)\\+(X_2\delta)(r-\lambda)\eta(X_1)-\delta^{2}[(X_{1}u)\eta(X_2)-(X_{2}u)\eta(X_1)]=0.
\end{split}
\end{equation}    
Now replacing $X_2$ by $\xi$ gives that\\
\begin{equation}\label{3.26}
    \delta^{2}[(X_{1}u)\eta(\xi)-(\xi u)\eta(X_1)]=0.
\end{equation}
Since $\delta\neq0$, so 
\begin{equation}\label{3.27}
    \nabla u=(\xi u)\xi.
\end{equation}
From equation \eqref{3.27} we can write $du=(\xi u)\eta$, where $d$ denotes the exterior differentiation. Taking the exterior derivative of the previous equation yields $d^{2} u=d(\xi u)\wedge\eta+(\xi u)d\eta$. Using \textit{Poincare lemma} $d^2\equiv 0$ and then executing wedge product with $\eta$ in the foregoing equation we obtain $(\xi u)\eta\wedge d\eta=0$.
Hence $\xi u=0$, since $\eta\wedge d\eta\neq 0$ in $M^{n}$. Thus we get $du=0$ which implies that $u$ is constant. This completes the proof.

\section{$\delta$-almost Yamabe soliton on a n-dimensional para-cosymplectic manifold.}
In 2004, Dacko \cite{Dac} introduced the notion of para-cosymplectic manifold. The
fundamental 2-form $\Phi$ is defined on an almost paracontact metric manifold $(M, \phi, \xi, \eta, g)$ by $\Phi(X_1,X_2)=g(X_1, \phi X_2)$ for any vector fields
$X_1$ and $X_2$ on $M$. It is evident that the skew-symmetricness of the 2-form $\Phi$ arises from $\phi$.\\
An almost paracontact metric manifold becomes almost para-cosympletic when the forms $\eta$ and $\Phi$ are closed, meaning $d\eta=0$ and $d\Phi=0$, respectively. Moreover, if the manifold also satisfies normality conditions, it is termed a para-cosymplectic manifold. Alternatively, we can describe an almost paracontact metric manifold is para-cosymplectic if the forms $\eta$ and $\Phi$ are parallel with respect to the Levi-Civita connection $\nabla$ of the metric, meaning $\nabla\eta=0$ and $\nabla\Phi=0$, respectively. Here are some important relations that hold true for any para-cosymplectic manifold.
\begin{align}
        R(X_1,X_2)\xi &= 0,\label{4.1}\\
        (\nabla_{X_1}\phi) &= 0,\label{4.2}\\
        \nabla_{X_1}\xi &= 0,\label{4.3}\\
        \mathcal{S}(X_1, \xi) &= 0,\label{4.4}\\
        Q\xi &= 0,\label{4.5}
\end{align}
where $X_1$ is an arbitrary vector field and $R, \nabla, S$ and $Q$ are the usual notations.\\
A vector field $Z$ is called a conformal Killing vector field, or simply a conformal vector field, if there exists a smooth function $\rho$ such that
\begin{equation}\label{4.6}
    \mathfrak{L}_{Z}g = 2\rho g.
\end{equation}
Here, $\rho$ is referred to as the conformal coefficient. When $\rho$ is set to zero, the conformal vector field reduces to a Killing vector field. 
\vspace{10pt}

Finally, we consider a n-dimensional para-cosymplectic manifold admitting $\delta$-almost Yamabe soliton.
\begin{theorem}
    If the metric $g$ of a n-dimensional para-cosymplectic manifold $(M, \phi, \xi, \eta, g)$ admits a conformal vector field $Z$, represents a $\delta$-almost Yamabe soliton then the scalar curvature of the manifold is $r=\lambda +\rho\delta$ and $\rho=(\mathfrak{L}_{Z}\eta)\xi =-\eta(\mathfrak{L}_{Z}\xi) $.
\end{theorem}
\proof Using the equation \eqref{4.6} and \eqref{1.3} for n-dimensional para-cosymplectic manifold we get
\begin{align}
     \delta\rho g(X_1,X_2)-(r-\lambda)g(X_1,X_2)&=0,\nonumber\\
     r=\lambda + \rho\delta .\label{4.10}
\end{align}
Now from the equation \eqref{3.37}, \eqref{3.38} and \eqref{4.10} , we obtain the following relation
\begin{equation}
    \rho=\frac{(r-\lambda)}{\delta}=(\mathfrak{L}_{Z}\eta)\xi =-\eta(\mathfrak{L}_{Z}\xi).\nonumber
\end{equation}
This completes the proof.

\begin{corollary}
    In n-dimensional para-cosymplectic manifold, a $\delta$-almost Yamabe soliton with conformal vector field is expanding, steady or shrinking according as, $r>\rho\delta$, $r=\rho\delta$ or $r<\rho\delta$.
\end{corollary}

\section{Examples}
Now, we provide some examples of $K$-paracontact metric manifold with $\delta$-almost Yamabe soliton.
\begin{exmp}
Let $M={(x_1,x_2,x_3)\in R^3}$, where $(x_1,x_2,x_3)$ are the standard coordinates in $R^3$, be a 3-dimensional manifold whose metric is $g$ with basis $\{e_1, e_2, e_3\}$ such that the Lie brackets are (see \cite{Kumara}, Example 14)
\begin{equation*}
[e_1,e_2]=-2e_3,\hspace{6pt} [e_3,e_2]=(u+1)e_1,\hspace{6pt} [e_3,e_1]=(u+1)e_2,
\end{equation*}
where $u$ is a smooth function.
We define a paracontact metric structure $(\phi,\xi,\eta,g)$ on $M$ by taking $\xi=e_3$ and 
\begin{equation*}
\phi e_1=e_2,\hspace{6pt} \phi e_2=e_1,\hspace{6pt} \phi \xi=0, 
\end{equation*}
\begin{equation*}
g(e_1,e_1)=g(\xi,\xi)=1,\hspace{6pt} g(e_2,e_2)=-1,\hspace{6pt}  \eta(X)=g(X,\xi),
\end{equation*}
for any $X\in\chi(M)$.
This paracontact metric structure is $K$-paracontact (in fact para-Sasakian) because $he_1=he_2=h\xi=0$.\\
Now by Koszul's formula, we get
\begin{equation*}
    \nabla_{e_1}e_1=0, \hspace{8pt}  \nabla_{e_1}e_2=-2\xi,\hspace{8pt} \nabla_{e_1}e_3=-(u+1)e_2,
\end{equation*}
\begin{equation}\label{5.1}
    \nabla_{e_2}e_1=2e_3, \hspace{8pt}\nabla_{e_2}e_2=0, \hspace{8pt} \nabla_{e_2}e_3=-(u+1)e_1,
\end{equation}

\begin{equation*}
    \nabla_{e_3}e_1=(u+1)e_2, \hspace{8pt}\nabla_{e_3}e_2=(u+1)e_1,\hspace{8pt}\nabla_{e_3}e_3=0.
\end{equation*}\\
It is implies that $M^3$ is $K$-paracontact manifold because $he_1=he_2=he_3=0$ and by a simple calculation gives scalar curvature $r=-2(2u+1)$.\\
Now if we consider $r$ is a nonzero constant scalar curvature (i.e., $u$ is constant) and $\delta=\lambda$ then from equation \eqref{1.3} we obtain 
\begin{equation}\label{5.2}
    g(\nabla_{X_1}Z,X_2)+g(X_1,\nabla_{X_2}Z)=2(\frac{r-\lambda}{\lambda})g(X_1,X_2).
\end{equation}
Replacing $Z$ by $\xi$  in \eqref{5.2} we get
\begin{equation*}
\begin{split}
    (\frac{r-\lambda}{\lambda})g(X_1,X_2)&= 0,\\
    \lambda &=r=constant.
\end{split}
\end{equation*}
It is easy to see that the metric $g$ satisfies the soliton equation. Therefore $M$ admits a $\delta$-Yamabe soliton.
\end{exmp}

\begin{exmp}
We define a normal almost paracontact metric structure $(\phi, \xi, \eta, g)$
on a Euclidean space $M = R^3$ with Cartesian coordinates $(x, y, z)$ as follows: (see \cite{Kumara}, Example 13)
\begin{equation*}
  \phi\Big(\frac{\partial}{\partial x}\Big)=-3z^2\frac{\partial}{\partial x}, \hspace{8pt}  \phi\Big(\frac{\partial}{\partial y}\Big)=3z^2\frac{\partial}{\partial y}, \hspace{8pt}  \phi\Big(\frac{\partial}{\partial z}\Big)=0, 
\end{equation*}
\begin{equation*}
\xi=\frac{\partial}{\partial z},\hspace{8pt} \eta=dz, \hspace{8pt} 
(g_{ij})=\Bigg(\begin{array}{ccc}
     exp(2z^3)&0&0 \\
     0&exp(-2z^3)&0\\
     0&0&1
     \end{array}\Bigg).
\end{equation*}
From Koszul's formula, the Levi-Civita connection $\nabla$ is given by
\begin{equation*}
    \nabla_{\partial x}\partial x=-3z^2exp(2z^3)\partial z, \hspace{8pt}  \nabla_{\partial x}\partial y=0,\hspace{8pt} \nabla_{\partial x}\partial z=3z^2\partial x,
\end{equation*}
\begin{equation}\label{5.3}
    \nabla_{\partial y}\partial x=0, \hspace{8pt}\nabla_{\partial y}\partial y=3z^2exp(-2z^3)\partial z, \hspace{8pt} \nabla_{\partial y}\partial z=-3z^2\partial y,
\end{equation}

\begin{equation*}
    \nabla_{\partial z}\partial x=3z^2\partial x, \hspace{8pt}\nabla_{\partial z}\partial y=-3z^2\partial y,\hspace{8pt}\nabla_{\partial z}\partial z=0,
\end{equation*}\\
where $\partial x=\frac{\partial}{\partial x}, \partial y=\frac{\partial}{\partial y}$ and $\partial z=\frac{\partial}{\partial z}$. This paracontact metric structure is
$K$-paracontact because $h(\partial x) = h(\partial y) = h(\partial z) = 0$ (being 3-dimensional it is para-Sasakian). Furthermore, a straightforward computation gives the scalar curvature $r = -18z^4$.\\
We suppose that $Z=(f_1+f_2)\partial x+(f_2+f_3)\partial y + (f_3+f_1)\partial z$ is an arbitrary potential vector field defined on $M$, where $f_1, f_2$ and $f_3$ are smooth functions. A $\delta$-almost Yamabe soliton equation \eqref{1.3} is written by
\begin{equation}\label{5.4}
\frac{\delta}{2}\{g(\nabla_{X}Z,Y ) + g(\nabla_{Y}Z,X)\} = (r -\lambda)g(X, Y ).
\end{equation}
As a result of \eqref{5.3} and \eqref{5.4}, we can conclude that $M$ admits a $\delta$-almost Yamabe soliton only if the smooth functions $f_1, f_2, f_3$ satisfy the following relations:
\begin{align*}
    \delta\Big[\frac{\partial(f_1+f_2)}{\partial x} + 3z^2(f_1+f_3)\Big] + 18z^4 + \lambda &= 0,\\
    \delta\Big[\frac{\partial(f_2+f_3)}{\partial y} - 3z^2(f_1+f_3)\Big] + 18z^4 + \lambda &= 0,\\
    \delta\Big[\frac{\partial(f_1+f_3)}{\partial z}\Big] + 18z^4 + \lambda &= 0,\\
    \frac{\partial(f_2+f_3)}{\partial x} exp(-2z^3)+ \frac{\partial(f_1+f_2)}{\partial y} exp(2z^3) &= 0,\\
    \frac{\partial(f_1+f_3)}{\partial y} + \frac{\partial(f_2+f_3)}{\partial z} exp(-2z^3) &= 0,\\
    \frac{\partial(f_1+f_2)}{\partial z} exp(2z^3)+ \frac{\partial(f_1+f_3)}{\partial x} &= 0.
\end{align*}
\end{exmp}
\section{Conflicts of interest}
On behalf of all authors, the corresponding author states that there is no conflict of interest.
\section{Data availability statement}
There is no data available for this submission.
\section{Authors contributions Statement}
Conceptualization, R. Biswas, S. Dey, A. Bhattacharyya; methodology, R. Biswas, S. Dey, A. Bhattacharyya;
writing review and editing, R. Biswas, S. Dey, A. Bhattacharyya. All authors have read and reviewed the
manuscript.
\section{Acknowledgements}
The first author is financially supported by UGC Junior Research Fellowship of India, NTA Ref. No.: 231610035035. The authors are very much thankful to the referee for his or her valuable comments and suggestions for improvements of this paper. 

\end{document}